\documentclass[12pt]{article}

\usepackage{fullpage,amsfonts,amsmath,amsthm,amssymb,mathrsfs,graphicx,epstopdf}
\usepackage[top=2.54cm, bottom=2.54cm, left=2.54 cm, right=2.54 cm]{geometry}
 \newcommand{\lab}[1]{\label{#1}}                

 \usepackage[usenames,dvipsnames]{color}

\newcommand{\remove}[1]{}
\newcommand\eqn[1]{(\ref{#1})}

\newcommand{\be}{\begin{equation}}
\newcommand{\bel}[1]{\begin{equation}\lab{#1}\ }
\newcommand{\ee}{\end{equation}}
\newcommand{\bea}{\begin{eqnarray}}
\newcommand{\eea}{\end{eqnarray}}
\newcommand{\bean}{\begin{eqnarray*}}
\newcommand{\eean}{\end{eqnarray*}}

\newtheorem{thm}{Theorem}
\newtheorem{cor}[thm]{Corollary}
\newtheorem{con}[thm]{Conjecture}
\newtheorem{lemma}[thm]{Lemma}

\newtheorem{claim}[thm]{Claim}

\newtheorem{remark}[thm]{Remark}

\def\proof{\noindent{\bf Proof.~}~}
\def\qed{~~\vrule height8pt width4pt depth0pt}

\def\B{{\mathcal B}}

\def\G{{\mathcal G}}

\def\N{{\mathcal N}}

\def\ex{{\mathbb E}}
\def\pr{{\mathbb P}}
\def\cov{{\bf Cov}}
\def\var{{\bf Var}}

\def\bfd{{\bf d}}

\def\bfg{{\bf g}}

\def\bfx{{\bf x}}

\def\eps{\epsilon}

\date{}

\title{The number of perfect matchings, and the nesting properties, of random regular graphs}
\author{Pu Gao \\
University of Waterloo\\
pu.gao@uwaterloo.ca 
}
\begin{document}
\maketitle

\begin{abstract}
We prove that the number of perfect matchings in $\G(n,d)$ is asymptotically normal when $n$ is even, $d\to\infty$ as $n\to\infty$, and $d=O(n^{1/7}/\log n)$. This is the first distributional result of spanning subgraphs of $\G(n,d)$ when $d\to\infty$.

Moreover, we prove that $\G(n,d-1)$ and $\G(n,d)$ can be coupled so that $\G(n,d-1)$ is a subgraph of $\G(n,d)$ with high probability when $d\to\infty$ and $d=o(n^{1/3})$. Further, if $d=\omega(\log^7 n)$, $d=O(n^{1/7}/\log n)$, and $d\le d'\le n-1$  then $\G(n,d)$ and $\G(n,d')$ can be coupled so that asymptotically almost surely (a.a.s.)  $\G(n,d)$ is a subgraph of $\G(n,d')$.
\end{abstract}

\section{Introduction}

In this paper we address two problems  regarding $\G(n,d)$, the random $d$-regular graph: the limiting distribution of the number of perfect matchings in $\G(n,d)$, and the sandwich conjectures of random regular graphs in terms of nesting $\G(n,d_1)\subseteq \G(n,d_2)$ with high probability when $d_1\le d_2$.
\subsection{The number of perfect matchings}
The study of subgraphs lies in the centre of random graph theory. The commonly studied examples include spanning subgraphs such as perfect matchings, Hamilton cycles, spanning trees, $H$-factors where $H$ has a fixed size, as well as smaller subgraphs such as independent sets, cycles, and in general subgraphs isomorphic to some given $H$ of fixed size.
Let $Z_H$ denote the number of subgraphs isomorphic to $H$. The phase transition of positive $Z_H$ and the distribution of $Z_H$ have been well studied in $\G(n,p)$ and in $\G(n,m)$ for both small and large $H$. It is interesting that $Z_H$ has different types of distributions for small and large $H$. If $H$ has fixed size and is balanced\footnote{A graph $G$ is balanced if $\max_{H\subseteq G}|E(H)|/|V(H)|=|E(G)|/|V(G)|$. Graph $G$ is said strictly balanced if $\max_{H\subseteq G, H\neq G}|E(H)|/|V(H)|<|E(G)|/|V(G)|$.}, then $Z_H$ is asymptotically normally distributed in $\G(n,p)$ and $\G(n,m)$ when $\ex Z_H\to\infty$~\cite{rucinski1988small}. The distribution of $Z_H$ for $H$ with linear size becomes complicated. For $p\gg n^{-1/2}$ and $p$ not too close to 1, the numbers of perfect matchings, Hamilton cycles, and spanning trees are log-normally distributed in $\G(n,p)$, but are normally distributed in $\G(n,m)$ when $m\gg n^{3/2}$~\cite{janson1994numbers}. For $m=\Theta(n^{3/2})$, these random variables also become log-normally distributed in $\G(n,m)$. It is not known if they remain log-normally distributed in $\G(n,p)$ when $p=O(n^{-1/2})$ and in $\G(n,m)$ when $m=o(n^{3/2})$, although it is conjectured so~\cite{janson1994numbers}. The distributional phase transition of $Z_H$ when the size of $H$ grows from constant to linear size has been studied in~\cite{gao2016transition}, where $H$ is the number of $\ell$-matchings (i.e.\ matchings of size $\ell$). Its distribution in $\G(n,p)$ changes from normal to log-normal at the critical value $\ell=\ell(p)\approx n\sqrt{p}$. Such distributional phase transition is also observed in $\G(n,d)$ when $d$ is a fixed constant. It is well known that the distributions of short cycles in $\G(n,d)$ are asymptotically Poisson~\cite{bollobas1980probabilistic,wormald1981asymptotic}, whereas the distribution of the number of large subgraphs such as perfect matchings and Hamilton cycles in $\G(n,d)$ is of an unusual type~\cite{janson1995random} as follows. Suppose that $Z$ is the number of perfect matchings (or the number of Hamilton cycles) in $\G(n,d)$. Then the limiting distribution of the logarithm of $Z/\ex Z$ is an infinite linear combination of independent Poisson variables. More precisely,
\begin{equation}
\frac{Z}{\ex Z} \to \prod_{i=1}^{\infty} (1+\delta_i)^{X_i} e^{-\lambda_i\delta_i}, \quad \mbox{as $n\to\infty$}, \lab{log-Poisson}
\end{equation}  
where $X_1,X_2,\ldots $ are independent Poisson variables with mean $\lambda_1,\lambda_2,\ldots$, and $\delta_1,\delta_2,\ldots$ are real numbers whose values depend on which subgraphs (i.e.\ perfect matchings or Hamilton cycles) $Z$ counts. The distribution of $Z$ is determined by using the small subgraph conditioning method, originally developed by Robinson and Wormald~\cite{robinson1992almost,robinson1994almost} to prove Hamiltonicity of $\G(n,d)$. The argument is then tuned to produce the distribution result of $Z$ by Janson~\cite{janson1995random}.  Recently, Greenhill, Isaev and Liang~\cite{greenhill2020spanning} proved that the number of spanning trees in $\G(n,d)$ has the same type of distribution as~\eqn{log-Poisson}. On the other hand, Garmo~\cite{garmo1999asymptotic} studied the distributional phase transition of the number of $\ell$-cycles in $\G(n,d)$ as $\ell$ grows from constant to linear in $n$. Its limiting distribution changes from a linear combination of independent Poisson variables to the exponential of that form, and the critical phase transition occurs when $\ell$ becomes linear in $n$.

It is natural to ask, in the case $d\to\infty$, whether the distribution type of these subgraphs (e.g.\ perfect matchings, Hamilton cycles, spanning trees) are the same as, or analogous to, that for constant $d$, and whether the distributional phase transitions occur when the size of the subgraphs (e.g. $\ell$-matchings and $\ell$-cycles) grows from constant to linear in $n$, as for constant $d$. We give a negative answer to this question.
There have been few distributional results that are known for the number of subgraphs of $\G(n,d)$ when $d\to\infty$, even for small subgraphs. The limiting distribution of the number of $\ell$-cycles was extended from constant $d$ and $\ell$ to those such that $(d-1)^{2\ell-1}=o(n)$ by McKay, Wormald and Wysocka~\cite{mckay2004short}.  Z. Gao and Wormald~\cite{gao2008distribution} determined the limiting distributions of strictly balanced graphs of fixed sizes for $d$ that grows sufficiently slowly with $n$. There has been no result on the distribution of the number of subgraphs whose size is beyond $\log n$ when $d\to\infty$. In particular, the analysis for the number of perfect matchings, Hamilton cycles, and the spanning trees when $d=O(1)$, based on the configuration model, cannot be extended easily to $d\to\infty$. 

One may expect that the number of large subgraphs such as perfect matchings or Hamilton cycles would be of log-normal type in $\G(n,d)$ as $d\to\infty$,  which can be viewed as an analog of~\eqn{log-Poisson}.  It is also reasonable to believe that the number of $\ell$-matchings may exhibit a distributional phase transition as $\ell$ grows from constant to linear in $n$, as that is what happens for the $\ell$-matchings in $\G(n,p)$ and for the $\ell$-cycles in $\G(n,d)$ for constant $d$. In contrast with the intuition, we show in this paper that the number of perfect matchings is asymptotically normally distributed in $\G(n,d)$ when $d\to\infty$ as $n\to\infty$ and $d=O( n^{1/7}/\log n)$. The power of the logarithmic term is not optimised.

\begin{thm}~\lab{thm:main}
Let $Y$ denote the number of perfect matchings in $\G(n,d)$ where $n$ is even. Then, $Y$ is asymptotically normally distributed if $d\to\infty$ as $n\to\infty$ and $d=O(n^{1/7}/\log n)$. More formally, 
\[
\frac{Y-\ex Y}{\sqrt{\var Y}} \xrightarrow{d} \N(0,1),\quad \mbox{as $n\to \infty$}.
\]
\end{thm}

This result suggests that there is likely no distributional phase transition on the number of $\ell$-matchings as $\ell$ grows. The condition 
$d=O( n^{1/7}/\log n)$ in the result is imposed only for technical reasons and we believe that the same distribution holds for all $d\to\infty$ until $d$ is too close to $n-1$; see Conjecture~\ref{con:normal} below.

To our knowledge, this is the first result on the limiting distribution of the number of spanning subgraphs in $\G(n,d)$ when $d\to\infty$. The main contribution of Theorem~\ref{thm:main} is the discovery of the distribution type of the number of perfect matchings, and we believe that this phenomenon is ubiquitous among other spanning subgraphs such as the number of Hamilton cycles. For future research, it would be interesting to determine the limiting distributions of the number of $\ell$-matchings and $\ell$-cycles in $\G(n,d)$ for all $\ell$. 

\begin{con}\label{con:normal}
The numbers of perfect matchings, Hamilton cycles,  spanning trees, and $k$-factors, where $k\le d-1$, are all asymptotically normally distributed in $\G(n,d)$ for all $d$ where $dn$ is even and $\min\{d,n-d\}\to\infty$ as $n\to\infty$ (and also $n$ is even in the case of perfect matchings, and $kn$ is even in the case of $k$-factors).
\end{con}

\begin{con}
Suppose $\min\{d,n-d\}\to\infty$ as $n\to\infty$. The number of $\ell$-cycles in $\G(n,d)$ is asymptotically normal for all $3\le \ell\le n$. The number of $\ell$-matchings of $\G(n,d)$ is asymptotically normal for all $3\le \ell\le n/2$.
\end{con}

\begin{remark}
The condition $n-d\to\infty$ in the above conjectures is likely not necessary. Indeed, when $n-d=o(n^{1/3})$ and $k=o(n^{1/3})$ the asymptotic number of $k$-factors of a $d$-regular graph $G$ is independent of $G$ and can be obtained by Theorem~\ref{thm:subgraph} in Section~\ref{sec:tools}.
\end{remark}

\subsection{The sandwich conjectures of $\G(n,d)$}

Analysis in $\G(n,d)$ is highly nontrivial, especially when $d\to\infty$. Kim and Vu initiated the study of approximating $\G(n,d)$ by $\G(n,p)$, known as the sandwich conjecture~\cite{kim2004sandwiching}. Since then several groups of authors~\cite{dudek2017embedding,gao2020sandwiching,klimovsova2020sandwiching,gao2020kim} have worked on this conjecture, and it is close to being fully resolved. Along the line of the research there has been new conjectures that are proposed, one of which is stated as follows~\cite[Conjecture 1.2]{gao2020sandwiching}.

\begin{con}\lab{con:sandwich}
Let $0\le d_1\le d_2\le n-1$ be integers, other than $(d_1,d_2)=(1,2)$ or $(d_1,d_2)=(n-3,n-2)$. Assume that $d_1n$ and $d_2n$ are both even. Then, there exists a coupling $(G_1,G_2)$ such that $G_1\sim \G(n,d_1)$, $G_2\sim \G(n,d_2)$ and $\pr(G_1\subseteq G_2)=1-o(1)$.
\end{con}
The conjecture is only known to be true for $(d_1,d_2)$ where $d_1=1$ and $3\le d_2\le n-1$, as well as for $(d_1,d_2)$ where $d_2-d_1$ is larger than some function of $d_1$ (see~\cite[Corollary 1.7]{gao2020sandwiching}  for the precise statement). When $d_1$ and $d_2$ are both fixed constants and $(d_1,d_2)\neq (1,2)$, it is known that $\G(n,d_2)$ is contiguous to the union of two independent copies of $\G(n,d_1)$ and $\G(n,d_2-d_1)$ conditional on $\G(n,d_1)$ and $\G(n,d_2-d_1)$ being disjoint (see more contiguity results in~\cite[Section 4]{wormald1999models}). However, contiguity does not imply a coupling as in the conjecture. 
In this paper we prove Conjecture~\ref{con:sandwich} for a certain range of $d_1$.
\begin{thm}\lab{thm:sandwich0}
Conjecture~\ref{con:sandwich} holds for all integers $d_1\le d_2\le n-1$ where $d_1=\omega(\log^7 n)$ and $d_1=O(n^{1/7}/\log n)$ if $n$ is even.
\end{thm}

Theorem~\ref{thm:sandwich0} follows as a corollary of~\cite[Theorem 2]{gao2020kim} and the following theorem that simultaneously couples a sequence of random regular graphs.

\begin{thm}\lab{thm:sandwich}
Suppose $d\to\infty$ and $d=O(n^{1/7}/\log n)$. For any $\eps_n=o(1)$, there is a multiple coupling $(G_{d}, G_{d+1},\ldots, G_{\lfloor(1+\eps_n)d\rfloor})$ such that marginally $G_i\sim \G(n,i)$ for all $d\le i\le \lfloor(1+\eps_n)d\rfloor$ and jointly $G_d\subseteq G_{d+1}\subseteq \cdots \subseteq G_{\lfloor(1+\eps_n)d\rfloor}$ a.a.s.. 
\end{thm}

If we only consider $(d_1,d_2)$ where $d_2=d_1+1$ then we have the following coupling theorem which holds for a much larger range of $d_1$.
\begin{thm}\lab{thm:sandwich2}
Suppose $d\to\infty$ and $d=o(n^{1/3})$. There is a coupling $(G_d,G_{d+1})$ where marginally $G_d\sim \G(n,d)$ and $G_{d+1}\sim \G(n,d+1)$, and  jointly $G_d\subseteq G_{d+1}$ a.a.s..
\end{thm}

Theorem~\ref{thm:sandwich2} follows as a corollary of a more general version (Theorem~\ref{thm:sandwich3}) which we state in Section~\ref{sec:sandwich}. Indeed, it is possible to prove that Theorem~\ref{thm:sandwich2} holds for all $d\to\infty$ and $d=o(n^{1/2})$. However, for the sake of a simpler proof, we did not pursue that. See Remark~\ref{remark} for more explanations. 

The two problems studied in this paper seem unrelated. However, the key ingredient in the proof of Theorem~\ref{thm:sandwich3} is the construction of a coupling procedure of $\G(n,d)$ and $\G(n,d+1)$. The success of the coupling relies on the concentration of the number of perfect matchings in $\G(n,d+1)$, which is one of the main results we obtain for the  first problem under study.

All asymptotics in the paper refers to $n\to\infty$.  Given two sequences of real numbers $a_n$ and $b_n$, we say $a_n=O(b_n)$ if there exists a constant $C>0$ such that $|a_n|\le C|b_n|$ for all $n$. We say $a_n=o(b_n)$ where $b_n>0$ for all sufficiently large $n$, if $\lim_{n\to\infty} a_n/b_n=0$. We say $a_n=\omega(b_n)$ if both $a_n$ and $b_n$ are positive for all sufficiently large $n$, and $b_n=o(a_n)$. We say $a_n=\Omega(b_n)$ if both $a_n$ and $b_n$ are positive for all sufficiently large $n$, and $b_n=O(a_n)$. 

\section{Proof of Theorem~\ref{thm:main}}

Recall that $Y$ denotes the number of perfect matchings in $\G(n,d)$. Throughout the paper we assume that $n$ is even.
 Let $X$ denote the number of triangles in $\G(n,d)$. We will approximate $Y$ by a linear function of $X$ using linear regression, and then study the distribution of $Y$ via analysing the distribution of $X$. This method is known as orthogonal decomposition and projection, developed by Janson~\cite{janson1994orthogonal}. Originally, it is developed to determine the limiting distribution of the number of (large) subgraphs in $\G(n,p)$, and Janson also applied the method to determine the distributions of the numbers of spanning trees, perfect matchings and Hamilton cycles in $\G(n,m)$~\cite{janson1994numbers}. We are not aware of any previous applications in $\G(n,d)$. More specifically, we approximate $Y$ by $Y^\star=aX+b$ where $a=\cov(X,Y)/\var X$ and $b=\ex Y-a\ex X$. The values of $a$ and $b$ are chosen so that $\ex (Y-Y^\star)=0$ and $\ex((Y-Y^\star)^2)$ is minimised. It follows immediately that $Y-Y^\star$ and $Y^\star$ are orthogonal random variables. We prove that $\ex((Y-Y^\star)^2)$ is sufficiently small and thus the distribution of $Y$ is asymptotically determined by the distribution of $aX+b$. Since $X$ is asymptotically normally distributed, so is $Y$. The expectation $\ex X$, the variance $\var X$ and the limiting distribution of $X$ have been studied in~\cite[Theorems 8 and 10]{gao2020triangles}, which we state below.

\begin{thm}\lab{thm:triangle}
Suppose $d=o(n^{2/5})$ and $d\ge 2$. Then,
\begin{align*}
&\ex X=\frac{(d-1)^3}{6}\left(1+O(1/n)\right),\quad\var X\sim \ex X,\\
& \frac{X-\ex X}{\sqrt{\var X}} \xrightarrow{d} \N(0,1) \ \mbox{as $n\to\infty$, provided $d\to\infty$}.
\end{align*}
\end{thm}

Next, we calculate the expectation $\ex Y$, the second moment $\ex Y^2$ and the covariance $\cov(X,Y)$, which allow us to estimate $a$ and $b$ and to bound $\ex((Y-Y^\star)^2)$.

\begin{thm}\lab{thm:2ndMoment} Suppose  $d=o(n^{1/3})$ and $d\ge 3$. Then,
\begin{align}
\ex Y&=\frac{n! }{(n/2)! 2^{n/2}} \left(\frac{e}{n}\right)^{n/2} \left(\frac{d-1}{d}\right)^{(\frac{d-1}{2})n}d^{\frac{n}{2}}\exp\left(\frac{1}{4}+O\left(\frac{d^3}{n}\right)\right) \label{Yexpectation}\\
\ex Y^2&=\left(1+\frac{1}{6d^3}+O\left(d^{-4}+\frac{d^3}{n}+\sqrt{\frac{d}{n}}\log^3n\right)\right) (\ex Y)^2 \label{Y2expectation}\\
\cov(X,Y)&=\left(-\frac{1}{d^3}+O\left(d^{-4}+\frac{d}{n}\right)\right)\ex X\ex Y.\label{XYcovariance}
\end{align}
\end{thm}

\begin{remark}
Although Theorem~\ref{thm:2ndMoment} holds for constant $d$, the expressions in~\eqn{Y2expectation} and~\eqn{XYcovariance} do not provide any asymptotic information about $\ex Y^2$ and $\cov(X,Y)$ since the error terms in $O(\cdot)$ are too big.
\end{remark}

From Theorem~\ref{thm:2ndMoment} it follows immediately that
\[
\var Y=\left(\frac{1}{6d^3}+O\left(d^{-4}+\frac{d^3}{n}+\sqrt{\frac{d}{n}}\log^3n\right)\right) (\ex Y)^2.
\]


\noindent {\em Proof of Theorem~\ref{thm:main}.\ } Recalling that $a=\cov(X,Y)/\var X$, we make the following claim.
\begin{claim}\lab{claim:Ydiff}
$\ex((Y-Y^\star)^2)=o(\var Y^\star)$.
\end{claim}
Since
\begin{eqnarray*}
Y&=& Y^\star +(Y-Y^\star), \quad \ex Y^\star = \ex Y,
\end{eqnarray*}
it follows then that
\begin{eqnarray}
\frac{Y-\ex Y}{\sqrt{\var Y^\star}}&=& \frac{Y^\star-\ex Y^\star}{\sqrt{\var Y^\star}} + \frac{Y-Y^\star}{\sqrt{\var Y^\star}} =  \frac{X-\ex X}{\sqrt{\var X}} + \frac{Y-Y^\star}{\sqrt{\var Y^\star}}.\label{Y-X-approx}
\end{eqnarray}
By Claim~\ref{claim:Ydiff}  and Markov's inequality, a.a.s.
\begin{equation}
|Y-Y^\star|=o(\sqrt{\var Y^\star}). \lab{small-error}
\end{equation}
Moreover, by the orthogonality of $Y-Y^\star$ and $Y^\star$, $\cov(Y-Y^\star,Y^\star)=0$ and so $\var Y=\var Y^\star+\var (Y-Y^\star) \sim \var Y^\star$ by Claim~\ref{claim:Ydiff}.
Thus, the left hand side of~\eqn{Y-X-approx} is asymptotic to $(Y-\ex Y)/\sqrt{\var Y}$ in probability, and the right hand side converges to a random variable whose limiting distribution is $\N(0,1)$ by~\eqn{small-error} and Theorem~\ref{thm:triangle}.
Consequently,
\[
\frac{Y-\ex Y}{\sqrt{\var Y}}\xrightarrow{d} \N(0,1),\quad \mbox{as $n\to\infty$}.\qed
\]

\noindent {\em Proof of Claim~\ref{claim:Ydiff}.\ } Since $\ex Y^\star=\ex Y$, we know that $\ex((Y-Y^\star)^2)=\var (Y-Y^\star)$. Moreover,  since $\cov(Y^\star, Y-Y^\star)=0$, we have $\var Y=\var Y^\star +\var (Y-Y^\star)$. Hence, it is sufficient to prove that
\begin{equation}
\var Y\sim \var Y^\star=\frac{\cov(X,Y)^2}{\var X}. \lab{asymp-equation1}
\end{equation}
By~\eqn{Yexpectation} and~\eqn{Y2expectation},
\[
\var Y=\ex Y^2-(\ex Y)^2 = \left(\frac{1}{6d^3}+O\left(d^{-4}+\frac{d^3}{n}+\sqrt{\frac{d}{n}}\log^3n\right)\right) (\ex Y)^2,
\]
and by~\eqn{XYcovariance} and Theorem~\ref{thm:triangle},
\[
\frac{\cov(X,Y)^2}{\var X}\sim \frac{\left(-\frac{1}{d^3}+O(d^{-4}+d/n)\right)^2(\ex X\ex Y)^2}{\ex X}\sim \frac{1}{6d^3}(\ex Y)^2. 
\]
Now~\eqn{asymp-equation1} follows since $d\to\infty$ and $d=O(n^{1/7}/\log n)$. \qed

\section{Proof of Theorem~\ref{thm:2ndMoment}}
\lab{sec:tools}

We will use the tools from~\cite[Theorem 1]{gao2020subgraph} and~\cite[Theorem 4.6]{mckay1985asymptotics} to
estimate $\ex Y$, $\ex Y^2$ and $\cov(X,Y)$.

\subsection{Edge and subgraph probabilities in $\G(n,d)$}

Let $H$ be a graph on $[n]$ and let $\bfd^H=(d_1^H,\ldots, d_n^H)$ denote the degree sequence of $H$. Suppose that $d_i^H\le d$ for every $1\le i\le n$.  Let $|H|$ denote the number of edges in $H$. The following result is a special case of~\cite[Theorem 1]{gao2020subgraph} for the conditional edge probability $\pr(uv\in\G(n,d)\mid H\subseteq \G(n,d))$.  
\begin{thm}\lab{thm:conditional} Suppose $d=o(n)$ and suppose that $H$ is a graph on $[n]$ such that $d_i^H\le d$ for every $1\le i\le n$ and $dn-2|H|=\Omega(dn)$. Then, 
\begin{eqnarray*}
&&\pr\Big(uv\in\G(n,d)\mid H\subseteq \G(n,d)\Big) =\left(1+O\left(\frac{d}{n}\right)\right)\frac{(d-d_u^H)(d-d_v^H)}{dn-2|H|}.
\end{eqnarray*}
\end{thm}
We will apply Theorem~\ref{thm:conditional} to estimate the probabilities of small subgraphs of $\G(n,d)$. For large subgraphs,
we apply instead the following enumeration result of McKay~\cite[Theorem 4.6]{mckay1985asymptotics}.
\begin{thm}\lab{thm:subgraph}
Let $\bfg=(g_1,\ldots,g_n)$ be a sequence of non-negative integers. Let $m=m(\bfg)=\|\bfg\|_1/2$. Let $X$ be a simple graph on $[n]$ with degree sequence $\bfx$. Let $\Delta(\bfg)$ and $\Delta(\bfx)$ denote the maximum components of $\bfg$ and $\bfx$ respectively. Suppose $\Delta(\bfg)\ge 1$, $\hat\Delta(\bfg)=o(m)$ where $\hat\Delta(\bfg)=\Delta(\bfg)^2+\Delta(\bfg)\Delta(\bfx)$. Define
\[
\lambda=\lambda(\bfg)=\frac{1}{4m(\bfg)}\sum_{j=1}^ng_i(g_i-1),\quad \mu=\mu(\bfg,X)=\frac{1}{2m(\bfg)}\sum_{ij\in X}g_ig_j.
\]
Let $N(\bfg,X)$ denote the number of simple graphs with degree sequence $\bfg$ and with no edge in common with $X$. Then, 
\[
N(\bfg, X)=\frac{(2m)!}{m!2^m \prod_{j=1}^n g_i!} \exp\Big(-\lambda(\bfg)-\lambda(\bfg)^2-\mu(\bfg,X)+O(\hat\Delta(\bfg)^2/m(\bfg))\Big).
\] 
\end{thm}

\begin{cor}\lab{cor:subgraph} Suppose $d\ge 3$ and $d=o(n)$.
Let $0\le k\le n/2$ be an integer.
Let $H$ be a graph on $[n]$ containing $k$ isolated edges and a collection of disjoint cycles spanning the remaining $n-2k$ vertices. Then, with $\alpha=2k/n$,
\begin{align*}
\pr(H\subseteq \G(n,d)) &= \frac{((d-2)n+2k)! \frac{dn}{2}! 2^{n-k} d^n(d-1)^{n-2k}}{(\frac{(d-2)n}{2}+k)!(dn)!} \exp\left(\phi(d,\alpha)+O\left(\frac{d^3}{n}\right)\right)\\
&= \left(\frac{e}{n}\right)^{(1-\frac{\alpha}{2})n} \left(\frac{d-2+\alpha}{d}\right)^{(\frac{d-2+\alpha}{2})n}d^{\frac{\alpha}{2}n}\left(d-1\right)^{(1-\alpha)n}\exp\left(\phi(d,\alpha)+O\left(\frac{d^3}{n}\right)\right),
\end{align*}
 where
 \begin{equation}
 \phi(d,\alpha)=\frac{4(d-2)^2-(d^2-5)\alpha^2-(2d^2-14d+20)\alpha}{4(d-2+\alpha)^2}.\lab{phi}
 \end{equation}
\end{cor}

\begin{remark}\lab{remark} (a) Theorem~\ref{thm:conditional} can be deduced from an earlier work than~\cite{gao2020subgraph}, e.g.\ by McKay~\cite{mckay1981subgraphs}. We cite~\cite[Theorem 1]{gao2020subgraph} as it is written in form of conditional edge probabilities, which is what we need in this paper.

(b) A stronger version of Theorem~\ref{thm:conditional} is available in~\cite[Theorem 6]{gao2020triangles} which estimates the conditional edge probabilities up to a relative error $d^2/n^2$ instead of $d/n$. Using that result, we can deduce Corollary~\ref{cor:subgraph}  with a smaller error $O(d^2/n)$ than $O(d^3/n)$. This will result in an improvement in the range of $d$ in several of theorems in the paper, e.g.\ in Theorems~\ref{thm:sandwich2},~\ref{thm:2ndMoment}, and~\ref{thm:sandwich3}. However, applying~\cite[Theorem 6]{gao2020triangles} involves more intensive calculations, and for a simpler proof we deduce Corollary~\ref{cor:subgraph} from Theorem~\ref{thm:subgraph} instead.

(c) It might be useful to notice, in applications of Corollary~\ref{cor:subgraph}, that $\phi(d,\alpha)$ is essentially $O(1)$. In particular, as $d\to\infty$, $\phi(d,\alpha)\to 1-(\alpha^2+2\alpha)/4$.
\end{remark}

{\em Proof of Corollary~\ref{cor:subgraph}.\ } Let $\bfd^H$ denote the degree sequence of $H$ and let $\bfg=\bfd-\bfd^H$ where $\bfd=(d,\ldots, d)$. Then, $\bfg$ has exactly $2k$ components of value $d-1$ and $n-2k$ components of value $d-2$. Hence, 
\begin{align*}
2m(\bfg)&=(d-2)n+2k\\
\lambda(\bfg)&=\frac{1}{2((d-2)n+2k)}\Big((d-1)(d-2)\cdot 2k+(d-2)(d-3)(n-2k)\Big)\\
\mu(\bfg,H)&=\frac{1}{(d-2)n+2k}\Big((d-1)^2\cdot k+ (d-2)^2(n-2k) \Big).
\end{align*}
and
\[
2m(\bfd)=dn,\quad \lambda(\bfd)=\frac{1}{2dn}\left(d(d-1)n\right),\quad \mu(\bfd,\emptyset)=0.
\]
Moreover,
\[
\hat\Delta(\bfg),\hat\Delta(\bfd)=O(d^2)\quad\mbox{and}\quad m(\bfg),m(\bfd)=\Omega(dn).
\]
Thus,
\begin{align}
\pr(H\subseteq\G(n,d))&=\frac{N(\bfg, H)}{N(\bfd,\emptyset)\nonumber}\\
&{\hspace{-2cm}}=\frac{((d-2)n+2k)!/\left(\frac{(d-2)n+2k}{2}! 2^{\frac{(d-2)n+2k}{2}}\right) }{(dn)!/\left(\frac{dn}{2}! 2^{\frac{dn}{2}}\right)}\cdot d^n (d-1)^{n-2k}\exp(\phi(d,\alpha)+O(d^3/n)), \lab{probability}
\end{align}
where
\begin{align}
\phi(d,\alpha)&=-\lambda(\bfg)-\lambda(\bfg)^2-\mu(\bfg,H)+\lambda(\bfd)+\lambda(\bfd)^2+\mu(\bfd,\emptyset)\nonumber\\
&=\frac{4(d-2)^2-(d^2-5)\alpha^2-(2d^2-14d+20)\alpha}{4(d-2+\alpha)^2}.\lab{phi-value}
\end{align}
See Maple calculations of~\eqn{phi-value} in the Appendix. Now the corollary follows by applying the Stirling formula to the factorials in~\eqn{probability}. The relative error $O(1/dn)$ in the Stirling formula is absorbed by $O(d^3/n)$.\qed
\subsection{$\cov(X,Y)$}

Fix a perfect matching $H$ of $K_n$. There are $(n/2)\cdot (n-2)$ ways to choose a triangle $T$ such that $|H\cap T|=1$. For any such $T$, by Theorem~\ref{thm:conditional} twice, the conditional probability of $T\subseteq \G(n,d)$ given $H\subseteq \G(n,d)$ is $(d-1)^2(d-1)_2/(M-n)^2 (1+O(d/n))$. 
There are $\binom{n}{3}-(n/2\cdot (n-2)) =(1+O(1/n)) n^3/6$ ways to choose a triangle $T$ such that $H\cap T=\emptyset$. For any such $T$, again by Theorem~\ref{thm:conditional} three times, the conditional probability of $T\subseteq \G(n,d)$ given $H\subseteq \G(n,d)$ is  $(d-1)_2^3/(M-n)^3 (1+O(d/n))$. 
Hence,
\begin{align*}
\ex XY&=\sum_{H\in \Phi} \pr(H\subseteq \G(n,d))\sum_{T\in\Psi}\pr(T\subseteq \G(n,d)\mid H\subseteq \G(n,d)),
\end{align*}
where $\Phi$ is the set of all perfect matchings in $K_n$, and $\Psi$ is the set of all triangles in $K_n$. By the discussions above, $\sum_{T\in\Psi}\pr(T\subseteq \G(n,d)\mid H\subseteq \G(n,d))$ is the same for every perfect matching $H$. Noting that $\sum_{H\in \Phi} \pr(H\subseteq \G(n,d))=\ex Y$, we have, by setting $x=1/d$,
\begin{align}
\ex XY&= \ex Y \left(\frac{n^2}{2} \frac{(d-1)^3(d-2)}{(M-n)^2}   +\frac{n^3}{6}\cdot  \frac{(d-1)^3(d-2)^3}{(M-n)^3}  \right) (1+O(d/n))\nonumber\\
&=\ex Y \left(\frac{d^2}{2}(1-x)(1-2x)  +\frac{d^3}{6} (1-2x)^3   \right)(1+O(d/n))\nonumber\\
&=(1-1/d^3+O(1/d^4+d/n)) \ex X\ex Y,\lab{Exy}
\end{align}
where the last equation above is obtained by taking the product of $\frac{d^2}{2}(1-x)(1-2x)  +\frac{d^3}{6} (1-2x)^3 $ and $(\ex X)^{-1}=(1+O(1/n))\frac{6}{d^3(1-x)^3}$ from Theorem~\ref{thm:triangle}, and then taking the Taylor expansion of the product at $x=0$. We include the Maple expansion formulae in the Appendix.

\subsection{$\ex Y$}
Let $H$ be a perfect matching of $K_n$. By Corollary~\ref{cor:subgraph} with $k=n/2$ (i.e.\ $\alpha=1$), we have $\phi(d,\alpha)=1/4$ and thus,
\[
\pr(H\subseteq\G(n,d))=(1+O(d^3/n))\rho_1(n,d).
\]
where
\begin{equation}
\rho_1(n,d)=\left(\frac{e}{n}\right)^{n/2} \left(\frac{d-1}{d}\right)^{(\frac{d-1}{2})n}d^{\frac{n}{2}}\exp\left(\frac{1}{4}\right). \lab{rho}
\end{equation}
Hence,
\begin{equation}
\ex Y=(1+O(d^3/n))\frac{n!}{(n/2)! 2^{n/2}} \rho_1(n,d). \lab{Y-expectation}
\end{equation}

\subsection{$\ex Y^2$}
Let $0\le k\le n/2$ be an integer. Fix two perfect matchings $H_1$ and $H_2$ of $K_n$ such that $|H_1\cap H_2|=k$. Let $\alpha=\alpha(k)=2k/n$. Then, by Corollary~\ref{cor:subgraph},
\[
\pr(H_1\cup H_2\subseteq \G(n,d)) = (1+O(d^3/n))\rho_2(n,d,\alpha)
\]
where
\begin{align}
\rho_2(n,d,\alpha)&=\frac{((d-2)n+2k)! \frac{dn}{2}! 2^{n-k} d^n(d-1)^{n-2k}}{(\frac{(d-2)n}{2}+k)!(dn)!} \exp\left(\phi(d,\alpha)\right)\lab{rho2}\\
&=\left(\frac{e}{n}\right)^{(1-\frac{\alpha}{2})n} \left(\frac{d-2+\alpha}{d}\right)^{(\frac{d-2+\alpha}{2})n}d^{\frac{\alpha}{2}n}\left(d-1\right)^{(1-\alpha)n}\exp\left(\phi(d,\alpha)+O(1/dn)\right),
\end{align}
with $\phi(d,\alpha)$ defined in~\eqn{phi}.
Next we compute the number of pairs $(H_1,H_2)$ of perfect matchings of $K_n$ such that $|H_1\cap H_2|=k$.

\begin{lemma}\lab{lem:k} Let $0\le k\le n/2-2$ be an integer.
The number of pairs $(H_1,H_2)$ of perfect matchings of $K_n$ such that $|H_1\cap H_2|=k$ is
\[
(1+O((n-2k)^{-1})) \frac{n!}{2^kk!\sqrt{e\pi (n-2k)/2}}.
\]
\end{lemma}

\proof The exponential generating function\footnote{We refer the readers to~\cite[Part A]{ac} for enumeration by generating functions.} for bi-coloured alternating cycles (i.e.\ edges along the cycle have alternating colours) of length at least 4 is
\[
F(z)=\sum_{n=2}^{\infty} \frac{(2n)!}{2\cdot 2n}\cdot 2 \cdot \frac{z^{2n}}{(2n)!} = -\frac{1}{2}(\log(1-z^2)+z^2).
\]
Thus,
the number of pairs $(H_1,H_2)$ of disjoint perfect matchings of $K_{2m}$ is 
\[
(2m)!\cdot [z^{2m}] e^{F(z)}= (2m)!\cdot [z^{2m}]\frac{e^{-z^2/2}}{\sqrt{1-z^2}}
\]
We know $\frac{e^{-z^2/2}}{\sqrt{1-z^2}}=\frac{e^{-1/2}}{\sqrt{1-z^2}}+O((1-z)^{1/2})$ by expanding $e^{-z^2/2}$ at $z=1$. Hence, by the transferring theorem~\cite[Theorems VI.3 and VI.4]{ac} and the binomial theorem,
$$
[z^{2m}]e^{F(z)} =[z^{2m}]\frac{e^{-1/2}}{\sqrt{1-z^2}} +[z^{2m}]O((1-z)^{1/2})  = \frac{e^{-1/2}}{4^m}\binom{2m}{m} +O(m^{-3/2}) = \frac{1}{\sqrt{e\pi m}}(1+O(m^{-1})).
$$
Thus, 
the number of pairs $(H_1,H_2)$ of perfect matchings of $K_n$ such that $|H_1\cap H_2|=k$ is
\[
\binom{n}{2k} \cdot \frac{(2k)!}{2^kk!} \cdot \frac{(n-2k)!}{\sqrt{e\pi (n-2k)/2}}(1+O((n-2k)^{-1}))=(1+O((n-2k)^{-1})) \frac{n!}{2^kk!\sqrt{e\pi (n-2k)/2}}.\qed
\]
By Lemma~\ref{lem:k} and recalling~\eqn{rho2},
\begin{align}
\ex Y^2&=\sum_{k=0}^{n/2-2} (1+O((n-2k)^{-1} +d^3/n)) \frac{n!}{2^kk!\sqrt{e\pi(n-2k)/2}} \rho_2(n,d,\alpha(k))+\ex Y. \lab{Ysquare}
\end{align}
Next, we show that the main contribution to $\ex Y^2$ comes from $k$ near some specific value. The proof of the lemma is postponed till Section~\ref{sec:sum-over-k}.
\begin{lemma}\lab{lem:sum-over-k} Assume $d=o(n^{1/3})$ and $d\ge 3$.
Let 
\[
\bar\alpha=\frac{1}{d},\quad \bar k=\lfloor \bar\alpha n/2\rfloor,\quad \bar\delta=\frac{2d}{n}\frac{d(d-2)}{(d-1)^2}.
\]
Then,
\[
\ex Y^2 = \left(1+O\left(\sqrt{\frac{d}{n}}\log^3n+\frac{d^3}{n}\right)\right)  \frac{2}{\sqrt{e\bar\delta}}  \frac{n!\rho_2(n,d,\bar\alpha)}{\bar k!2^{\bar k}\sqrt{n-2\bar k}}.
\]
\end{lemma}

\subsubsection{Comparing $\ex Y^2$ with $(\ex Y)^2$}
We complete the proof of Theorem~\ref{thm:2ndMoment} by verifying that
\begin{equation}
\ex Y^2=\left(1+\frac{1}{6d^3}+O(\xi)\right) (\ex Y)^2,\lab{relation}
\end{equation}
where $\xi=d^{-4}+\frac{d^3}{n}+\sqrt{\frac{d}{n}}\log^3 n$.
By~\eqn{Y-expectation} and Lemma~\ref{lem:sum-over-k},
\[
(\ex Y)^2=(1+O(\xi))2\left(\frac{n}{e}\right)^n\rho_1(n,d)^2,
\]
and
\[
\ex Y^2=(1+O(\xi)) \frac{2}{\sqrt{e\bar\delta}} \sqrt{\frac{n}{\bar k}}  \frac{(n/e)^n\rho_2(n,d,\bar\alpha)}{(2\bar k/e)^{\bar k}\sqrt{n-2\bar k}}.
\]
Hence,
\begin{align*}
\frac{\ex Y^2}{(\ex Y)^2}=(1+O(\xi)) \sqrt{\frac{n}{e\bar\delta\bar k(n-2\bar k)}}  \frac{\rho_2(n,d,\bar\alpha)}{(2\bar k/e)^{\bar k}\rho_1(n,d)^2}.
\end{align*}
By straightforward but tedious calculations (see Appendix for more details)
\begin{equation}
\sqrt{\frac{n}{e\bar\delta\bar k(n-2\bar k)}}  \frac{\rho_2(n,d,\bar\alpha)}{(2\bar k/e)^{\bar k}\rho_1(n,d)^2}=1+\frac{1}{6d^3}+O(\xi), \lab{asymp-equation}
\end{equation}
and now~\eqn{relation} follows.\qed

\subsubsection{Proof of Lemma~\ref{lem:sum-over-k}} \lab{sec:sum-over-k}
Recall from~\eqn{Ysquare} that
\begin{align}
\ex Y^2&=\sum_{k=0}^{n/2-2} (1+O((n-2k)^{-1} +d^3/n)) \frac{n!}{2^kk!\sqrt{e\pi(n-2k)/2}} \rho_2(n,d,\alpha(k))+\ex Y\nonumber\\
&=\sum_{k=0}^{n/2-2} (1+O((n-2k)^{-1} +d^3/n))\sqrt{\frac{2}{e\pi}} n! \varphi(k)+\ex Y, \lab{Ysquare2}
\end{align}
where
\[
\varphi(k)= \frac{\rho_2(n,d,\alpha(k))}{2^kk!\sqrt{n-2k}}.
\]
The proof of Lemma~\ref{lem:sum-over-k} is standard. We prove that the summand in~\eqn{Ysquare2} is maximised at $\bar k$. Then, we approximate the summation around $\bar k$ by an integral of a function of the form $e^{-x^2}$. The contributions to~\eqn{Ysquare2} from $k$ far away from $\bar k$ is negligible. 

It is easy to see then that $\exp(\phi(d,\alpha(k))-\phi(d,\alpha(k-1)))=\exp(O(1/n))$. Hence, by~\eqn{rho2},
\begin{equation}
\frac{\varphi(k)}{\varphi(k-1)}=\left\{
\begin{array}{ll}
\frac{(d-2)n+2k}{2(d-1)^2k}\left(1+O\left(\frac{1}{n}\right)\right)& \mbox{for all $k\le n/3$}\\
O\left(\frac{(d-2)n+2k}{2(d-1)^2k}\right)=O(1/d) &  \mbox{for all $n/3<k\le n/2-2$}
\end{array}
\right.\lab{ratio}
\end{equation}
By equating $\frac{(d-2)n+2k}{2(d-1)^2k}$ to 1 we find that $k=n/2d$. It follows immediately that
 at $\bar k$, the ratio $\frac{(d-2)n+2k}{2(d-1)^2k}$ is $1+O(d/n)$, since rounding $n/2d$ to $\bar k$ would change this ratio by $O(d/n)$. Moreover,  noticing that
\begin{eqnarray*}
\frac{(d-2)n+2\bar k+2i}{2(\bar k+i)(d-1)^2}&=&\frac{(d-2)n+2\bar k}{2\bar k(d-1)^2} \left(1+\frac{2i}{(d-2)n+2\bar k}\right)\left(1+\frac{i}{\bar k}\right)^{-1}\\
&=&\left(1+O(d/n)\right) \exp\left(\frac{2i}{(d-2)n+2\bar k}-\frac{i}{\bar k} +O\left(\frac{i^2}{\bar k^2}\right)\right),
\end{eqnarray*}
we have that for every positive $j=o(\bar k)$:
\begin{align}
\frac{\varphi(\bar k+j)}{\varphi(\bar k)}&=\left(1+O\left(\frac{j}{n}\right)\right)\prod_{i=1}^j \frac{(d-2)n+2\bar k+2i}{2(\bar k+i)(d-1)^2}\nonumber\\
&= \left(1+O\left(\frac{j}{n}\right)\right)\prod_{i=1}^j  \left(1+i\left(\frac{2}{(d-2)n+2\bar k}-\frac{1}{\bar k}\right) +O\left(\frac{i^2}{\bar k^2}+\frac{d}{n}\right)\right)\nonumber\\
&=\exp\left(-\bar\delta \sum_{i=1}^j i +O\left(\frac{j^3}{\bar k^2}+\frac{dj}{n}\right)\right)\nonumber\\
&=\exp\left(-\frac{\bar\delta}{2} j^2 +O\left(\bar\delta j+\frac{j^3}{\bar k^2}+\frac{dj}{n}\right)\right),\lab{ratio-error}
\end{align}
recalling that $\bar \alpha=1/d$, $\bar k=\lfloor \bar\alpha \cdot n/2\rfloor$ and
\begin{align*}
\bar\delta&=\frac{2d}{n}\frac{d(d-2)}{(d-1)^2}=\left(\frac{1}{\bar k}-\frac{2}{(d-2)n+2\bar k}\right) \left(1+O\left(\frac{1}{\bar k}\right)\right)\\
&=\left(\frac{1}{\bar k}-\frac{2}{(d-2)n+2\bar k}\right) \left(1+O\left(\frac{d}{n}\right)\right).
\end{align*}
Symmetric calculations show that~\eqn{ratio-error} holds also for every negative $j=o(\bar k)$.
It follows then that
\begin{align*}
\sum_{\bar k-\bar\delta^{-1/2}\log(1/\bar\delta) \le k\le \bar k+\bar\delta^{-1/2}\log(1/\bar\delta)} \varphi(k) 
&=\left(1+O\left(\xi \right)\right) \varphi(\bar k) \sqrt{\frac{2}{\bar\delta}} \int_{-\infty}^{+\infty} e^{-x^2} dx \\
&=\left(1+O\left(\xi \right)\right) \sqrt{\frac{2\pi}{\bar\delta}}\varphi(\bar k),
\end{align*}
where
\[
\xi=\bar\delta^{1/2} \log(1/\bar\delta) +\frac{ \bar\delta^{-3/2}\log^3(1/\bar\delta)}{\bar k^2} +\frac{d}{n} \bar \delta^{-1/2}\log(1/\bar\delta)+\bar\delta^{1/2}=O\left(\sqrt{\frac{d}{n}}\log^3n\right).
\]
Note that the first three terms in $\xi$ come from the accumulative error $O(\bar\delta j+j^3/\bar k^2+dj/n)$ in~\eqn{ratio-error}, and the last term comes from approximating the sum of $\exp(-\bar\delta j^2/2)$ by an integral.
The contributions to~\eqn{Ysquare2} from $k$ where $|k-\bar k|>\bar\delta^{-1/2}\log(1/\bar\delta)$ is smaller than $n^{-1}$ --- indeed, at most $\exp(-\Omega(\log^2 n))$ --- as a relative error by~\eqn{ratio} and standard calculations on summing a geometrically bounded series. 
So
\[
\ex Y^2= \left(1+O\left( \sqrt{\frac{d}{n}}\log^3n+\frac{d^3}{n}\right)\right)  \frac{2}{\sqrt{e\bar\delta}}  \frac{n!\rho_2(n,d,\bar\alpha)}{\bar k!2^{\bar k}\sqrt{n-2\bar k}},
\]
where the error $d^3/n$ is carried from~\eqn{Ysquare2}. \qed
\section{Proofs of Theorems~\ref{thm:sandwich0} -- \ref{thm:sandwich2}}
\lab{sec:sandwich}

\subsection{Proof of Theorem~\ref{thm:sandwich0}}
We prove Theorem~\ref{thm:sandwich0} assuming Theorem~\ref{thm:sandwich}. Given $0\le p\le n$, recall that $\G(n,p)$ denotes the Erd\H{o}s-R\'{e}nyi random graph with edge probability $p$. Let $d_1=\omega(\log^7 n)$. By~\cite[Theorem~2]{gao2020kim}, there exists $\delta_n=o(1)$ such that $\G(n,d_1)$, $\G(n, p_1)$ and $\G(n,d_2)$ can be coupled together, where $p_1=(1+\delta_n)d_1/n$ and $d_2\ge (1+2\delta_n)d_1$, such that a.a.s.\ $\G(n,d_1)\subseteq \G(n, p_1)\subseteq \G(n,d_2)$.  It follows now that Conjecture~\ref{con:sandwich} holds for any $(d_1,d_2)$ where $d_1=\omega(\log^7 n)$ and $d_2-d_1\ge 2\delta_n d_1$. Now suppose $d_1+1\le d_2<(1+2\delta_n)d_1$ and $d_1=O(n^{1/7}/\log n)$. Then, there is a coupling where a.a.s.\ $\G(n,d_1)\subseteq \G(n,d_2)$ by Theorem~\ref{thm:sandwich}. Now Theorem~\ref{thm:sandwich0} follows.\qed

\subsection{Couple $\G(n,d)$ and $\G(n,d+1)$}
Throughout this section we assume $d\to\infty$ and $d=o(n^{1/3})$. Our goal is to couple $\G(n,d)$ with $\G(n,d+1)$ so that $\G(n,d)\subseteq \G(n,d+1)$ with sufficiently high probability. In the next section, we ``stitch'' a sequence of such couplings together to obtain a simultaneous coupling as in Theorem~\ref{thm:sandwich}.

Given $\alpha=\alpha_n=o(1)$, define $\eta=\eta(\alpha)$ where
\[
\eta(\alpha)=2\alpha+\frac{1}{d^3\alpha^2}+\frac{C'd^3}{n\alpha^2}+\frac{C'\sqrt{d/n}\log^3n}{\alpha^2}, \quad\mbox{where $C'>0$ is a sufficiently large constant.}
\]
We prove the following stronger version of Theorem~\ref{thm:sandwich2}.
\begin{thm}\lab{thm:sandwich3} Assume $\alpha=o(1)$ is such that $\eta(\alpha)=o(1)$.
There is a coupling $(G_d,G_{d+1})$ where marginally $G_d\sim \G(n,d)$ and $G_{d+1}\sim \G(n,d+1)$, and jointly $G_d\subseteq G_{d+1}$ with probability at least $1-5\eta$ for all sufficiently large $n$. 
\end{thm}

\subsubsection{The coupling procedure}

In this section we assume that $\eta(\alpha)=o(1)$ for some $\alpha=o(1)$. This assumption ensures that various quantities in our procedure below are positive. This assumption also immediately implies that $d=o(n^{1/3})$ and $d\to\infty$, which satisfies the condition for Theorem~\ref{thm:2ndMoment}. In this subsection, with slight abuse of notation, we also let $\G(n,d)$ denote the set of $d$-regular graphs on $[n]$. From the context it is always clear whether we refer to a set of graphs or a random graph from the set.

For $G\in \G(n,d+1)$ let $Y(G)$ denote the number of perfect matchings of $G$. For $G\in \G(n,d)$ let $Z(G)$ denote the number of perfect matchings in $K_n\setminus G$. We say $G$ and $G'$ are related, denoted by $G\sim G'$, for $G\in \G(n,d)$ and $G'\in \G(n,d+1)$, if $G\subseteq G'$. We can represent this relation using an auxiliary directed bipartite graph ${\cal X}$ where $V({\cal X})=\G(n,d)\cup\G(n,d+1)$, and $(G,G')$ is an arc if $G\subseteq G'$. Thus, a $d$-regular graph $G$ in ${\cal X}$ has out-degree $Z(G)$, and a $(d+1)$-regular graph $G'$ in ${\cal X}$ has in-degree $Y(G')$. 

Let $\alpha=\alpha_n=o(1)$. By Theorem~\ref{thm:2ndMoment} and Chebyshev's inequality,
\begin{equation}
\pr_{\G(n,d+1)}\left(|Y-\ex Y|\ge  \alpha \ex Y \right) \le \frac{\var Y}{\alpha^2 (\ex Y)^2} =\frac{1}{6d^3\alpha^2}+O\left(\frac{1}{d^4\alpha^2}+\frac{d^3}{n\alpha^2}+\frac{\sqrt{d/n}\log^3n}{\alpha^2}\right). \lab{Y-conc}
\end{equation}
By Theorem~\ref{thm:subgraph} (with $\bfg$ being the all one vector and $X$ being a $d$-regular graph) and Theorem~\ref{thm:2ndMoment}, there exists a sufficiently large constant $C>0$ such that
\begin{align}
|Z(G)-Z^\star|&\le C\frac{d^2}{n}\cdot Z^{\star},\ \mbox{for all $d$-regular graph $G$},\lab{Z-conc}\\
|\ex Y(G)-Y^\star| &\le C\frac{d^3}{n}\cdot Y^{\star},\ \mbox{for $G\sim \G(n,d+1)$}, \lab{exY}
\end{align}
where
\[
Z^\star=\frac{n!e^{-d/2}}{(n/2)! 2^{n/2}},\quad Y^\star= \frac{n! e^{1/4} }{(n/2)! 2^{n/2}} \left(\frac{e}{n}\right)^{n/2} \left(\frac{d}{d+1}\right)^{\frac{dn}{2}}(d+1)^{\frac{n}{2}}.
\]
Let
\begin{align}
\underline{Y}&= (1-\alpha-Cd^3/n)\, Y^\star, \lab{Ybar} \\
\overline{Z}&= \left(1+C \frac{d^2}{n}\right)Z^\star. \lab{Zbar} 
\end{align}
Define
\begin{align}
\B&=\{G\in \G(n,d+1):\ Y(G)<\underline{Y} \} \lab{B}\\
\B'&=\{G\in \G(n,d+1): \ Y(G)>(1+\alpha+Cd^3/n)\, Y^\star\}.\lab{B'}
\end{align}
Since $\alpha=o(1)$ and $\eta(\alpha)=o(1)$, it follows immediately that $\alpha\gg d^3/n$. Thus, by~\eqn{Y-conc} and~\eqn{exY}, 
\begin{equation}
\pr(\B\cup\B')\le\frac{1}{6d^3\alpha^2}+O\left(\frac{1}{d^4\alpha^2}+\frac{d^3}{n\alpha^2}+\frac{\sqrt{d/n}\log^3n}{\alpha^2}\right).\lab{B-size}
\end{equation}
Let $D$ and $\hat D$ denote the total in-degrees of $\G(n,d+1)$ and $\G(n,d+1)\setminus \B$ respectively in ${\cal X}$. That is,
\begin{align}
D&=|\{(G,G')\in \G(n,d)\times\G(n,d+1): G\sim G'\}|\lab{D}\\
\hat D&=|\{(G,G')\in \G(n,d)\times(\G(n,d+1)\setminus \B): G\sim G'\}|. \lab{Dhat}
\end{align}
Further, let $d^-(\B)$ and $d^-(\B')$ denote the total in-degrees of $\B$ and $\B'$ respectively in ${\cal X}$. We prove the following bounds on  $d^-(\B)$ and $d^-(\B')$.
Recall that
\[
\eta(\alpha)=2\alpha+\frac{1}{d^3\alpha^2}+\frac{C'd^3}{n\alpha^2}+\frac{C'\sqrt{d/n}\log^3n}{\alpha^2},
\]
where $C'>0$ is a sufficiently large constant.
\begin{lemma}\lab{lem:degree-B} Assume $\alpha=o(1)$ is such that $\eta(\alpha)=o(1)$. Then,
$d^-(\B)+d^-(\B')\le \eta D$.
\end{lemma}

\proof 
The number of edges in ${\cal X}$ is
$
D=|\G(n,d+1)|\ex Y.
$
This can be rewritten as 
\[
d^-(\B)+d^-(\B')+|\G(n,d+1)\setminus(\B\cup \B')|\cdot \ex Y\left(1+\xi\right),
\]
where $|\xi|\le \alpha+O(d^3/n)$ since $|Y(G)/\ex Y-1|\le \alpha+O(d^3/n)$ for all $G\notin \B\cup\B'$ by the definition of $\B$ and $\B'$.
By~\eqn{B-size}, the above is equal to
\[
d^-(\B)+d^-(\B')+|\G(n,d+1)|\left(1-\xi'\right)\cdot \ex Y\left(1+\xi\right)
\]
where 
\[
0\le\xi'\le \frac{1}{6d^3\alpha^2}+O\left(\frac{1}{d^4\alpha^2}+\frac{d^3}{n\alpha^2}+\frac{\sqrt{d/n}\log^3n}{\alpha^2}\right).
\]
Thus,
\[
|\G(n,d+1)|\ex Y=d^-(\B)+d^-(\B')+|\G(n,d+1)|\cdot\ex Y\,(1-\xi'+\xi-\xi'\xi),
\]
which implies that
\[
d^-(\B)+d^-(\B')=|\G(n,d+1)|\cdot\ex Y(\xi'-\xi+\xi'\xi) <\eta D,
\]
where the last inequality holds because $\xi'-\xi+\xi'\xi < |\xi|+2\xi'<\eta$ by the definition of $\eta$. \qed
\medskip

Finally we are ready to define the coupling $(G_d,G_{d+1})$.

\begin{itemize}
\item Let $G_{d}$ be a uniformly random graph in $\G(n,d)$ and let $\bar G$ be the graph obtained from $G_{d}$ by adding a uniformly random perfect matching of $K_n\setminus G_{d}$. Let $H$ be a uniformly random graph in $\G(n,d+1)$ independent of $G_d$ and $\bar G$. 
\item If $\bar G\in \B$ then let $G_{d+1}=\bar G$ with probability $(1-\eta)\frac{Z(G_d)}{\overline{Z}} $ and let $G_{d+1}=H$ with the remaining probability.
\item If $\bar G\in \G(n,d+1)\setminus \B$, then 
\[
G_{d+1}=\left\{
\begin{array}{ll}
 \bar G & \mbox{with probability $(1-\eta)\frac{Z(G_d)}{\overline{Z}} \cdot \frac{\underline{Y}}{Y(G')}$}\\
 G'' & \mbox{with probability $(1-\eta)\frac{Z(G_d)}{\overline{Z}}  \frac{(\underline{Y}-Y(G''))}{\hat D}$ for every $G''\in\B$}\\
  H & \mbox{with the remaining probability}.
 \end{array}
 \right.
 \]
\end{itemize}
Note that $Z(G)\le \overline{Z}$ for every $G\in\G(n,d)$ by the definition of $\overline{Z}$. 
The following lemma bounds a few quantities in the above probability terms, and in particular, it justifies  that the coupling procedure is well defined (for all sufficiently large $n$).
\begin{lemma} \lab{lem:bounds} Assume $\alpha=o(1)$ is such that $\eta(\alpha)=o(1)$.
\begin{align*}
\frac{Z(G)}{\overline{Z}} &\ge 1-3C\frac{d^2}{n}, \quad\mbox{for every $d$-regular graph $G$,}\\  
\frac{\underline{Y}}{Y(G')}&\ge 1-3\alpha-3C\frac{d^3}{n}, \quad\mbox{for every $(d+1)$-regular graph $G'\in \G(n,d+1)\setminus(\B\cup\B')$,}
\end{align*}
and for every $(G,G')\in \G(n,d)\times \big(\G(n,d+1)\setminus\B\big)$ where $G\sim G'$,
\[
(1-\eta)\frac{Z(G)}{\overline{Z}} \cdot \frac{\underline{Y}}{Y(G')} + \sum_{G''\in\B} (1-\eta)\frac{Z(G)}{\overline{Z}}  \frac{(\underline{Y}-Y(G''))}{\hat D}\le 1.
\]
\end{lemma}

\proof Since $\alpha=o(1)$ and $\eta(\alpha)=o(1)$, it follows immediately that
$d=o(n^{1/3})$ and $d\to\infty$. The first inequality in the lemma follows by~\eqn{Z-conc} and~\eqn{Zbar}. The second inequality in the lemma follows by~\eqn{Ybar},~\eqn{B} and~\eqn{B'}.
For the last inequality, note that 
\[
\frac{Z(G)}{\overline{Z}}, \frac{\underline{Y}}{Y(G')} \le 1
\]
always by~\eqn{Z-conc} and the definition of $\B$. Thus it is sufficient to show that
\[
 \sum_{G''\in\B}  \frac{(\underline{Y}-Y(G''))}{\hat D}\le \eta.
\]
By~\eqn{Y-conc}, 
\[
|\B|\le\left(\frac{1}{6d^3\alpha^2}+O\left(\frac{1}{d^4\alpha^2}+\frac{d^3}{n\alpha^2}+\frac{\sqrt{d/n}\log^3n}{\alpha^2}\right)\right) |\G(n,d+1)|.
\]
 By Lemma~\ref{lem:degree-B}, 
\[
\hat D=(1+O(\eta))D= (1+O(\eta)) |\G(n,d+1)| \ex Y= (1+O(\eta)) |\G(n,d+1)| \underline{Y}.
\]  
Thus,
\[
 \sum_{G''\in\B}  (\underline{Y}-Y(G'')) \le \underline{Y} |\B| \le \left(\frac{1}{6d^3\alpha^2}+O\left(\frac{1}{d^4\alpha^2}+\frac{d^3}{n\alpha^2}+\frac{\sqrt{d/n}\log^3n}{\alpha^2}\right)\right)\cdot |\G(n,d+1)| \underline{Y} \le \eta\hat D,
\]
by the definition of $\eta$. Thus, the last inequality of the lemma follows. \qed

\subsubsection{Proof of Theorem~\ref{thm:sandwich3}}
 By the construction, $G_{d}$ is obviously distributed as $\G(n,d)$ marginally. We prove that the marginal distribution of $G_{d+1}$ is $\G(n,d+1)$.  
Let
\[
\sigma_d=\frac{1}{|\G(n,d)|},\quad \sigma_{d+1}=\frac{1}{|\G(n,d+1)|}.
\]
Let $\hat G$ be a $(d+1)$-regular graph.
If $\hat G\in\G(n,d+1)\setminus\B$ then
\begin{equation}
\pr(G_{d+1}=\hat G)=\sum_{G: G\sim \hat G} \frac{\sigma_d}{Z(G)} \cdot (1-\eta)\frac{Z(G)}{\overline{Z}} \cdot \frac{\underline{Y}}{Y(\hat G)} +\varphi, \lab{prob1}
\end{equation}
where 
\begin{align*}
\varphi=&\sum_{\substack{(G,G'):\\ G\sim G',\ G'\in \B}} \frac{\sigma_d}{Z(G)}\left(1-(1-\eta)\frac{Z(G)}{\overline{Z}} \right) \sigma_{d+1}\\
&+\sum_{\substack{(G,G'):\\ G\sim G',\ G'\notin \B}}\frac{\sigma_d}{Z(G)} \left(1-(1-\eta)\frac{Z(G)}{\overline{Z}} \cdot \frac{\underline{Y}}{Y(G')}- \sum_{G''\in\B}(1-\eta)\frac{Z(G)}{\overline{Z}}  \frac{(\underline{Y}-Y(G''))}{\hat D}\right)\sigma_{d+1}.
\end{align*}
In the first summation in~\eqn{prob1}, $\sigma_d/Z(G)$ is the probability that $G_d=G$ and $\bar G=\hat G$. Conditioning on that, $(1-\eta)\frac{Z(G)}{\overline{Z}} \cdot \frac{\underline{Y}}{Y(G')}$ is the probability that $G_{d+1}$ is set to be $\bar G$ . Thus this summation gives the contribution to $\pr(G_{d+1}=\hat G)$ from the case that $\bar G=\hat G$ and $G_{d+1}$ is set to be $\bar G$. Similarly, it is easy to see that
 $\varphi$ is the probability that $H=\hat G$ and $G_{d+1}$ is set to be $H$. 
 Note that the value of $\varphi$ is independent of $\hat G$.
 Hence, by noting that $Y(\hat G)=|\{G: G\sim \hat G\}|$, we obtain
 \[
 \pr(G_{d+1}=\hat G)= (1-\eta) \sigma_d \frac{\underline{Y}}{\overline{Z}}  +\varphi,
 \]
 which is independent of $\hat G$ for all $\hat G\in \G(n,d+1)\setminus\B$. 
 
 Next, suppose $\hat G\in \B$. Then,
 \[
 \pr(G_{d+1}=\hat G)=\sum_{G: G\sim \hat G} \frac{\sigma_d}{Z(G)} \cdot (1-\eta)\frac{Z(G)}{\overline{Z}} +\sum_{\substack{(G,G'):\\ G\sim G',\ G'\notin \B}}  \frac{\sigma_d}{Z(G)}\cdot (1-\eta)\frac{Z(G)}{\overline{Z}}  \frac{(\underline{Y}-Y(\hat G))}{\hat D} +\varphi,
 \]
where the second summation above is from the case where $G_d=G$, $\bar G=G'\notin \B$,  and $G_{d+1}$ is set to be $G''=\hat G$ which occurs with probability  $(1-\eta)\frac{Z(G)}{\overline{Z}}  \frac{(\underline{Y}-Y(\hat G))}{\hat D}$, given $(G,G')$.
 The first summation above gives $(1-\eta)\sigma_d Y(\hat G)/\overline{Z}$. The second summation above gives
 $(1-\eta)\sigma_d (\underline{Y}-Y(\hat G))/\overline{Z}$ by~\eqn{Dhat}. Hence,
 \[
 \pr(G_{d+1}=\hat G)= (1-\eta) \sigma_d \frac{\underline{Y}}{\overline{Z}}  +\varphi,
 \]
   which is independent of $\hat G$ for all $\hat G\in \B$, and is the same for all $\hat G\in \G(n,d+1)\setminus \B$. This confirms that the marginal distribution of $G_{d+1}$ is uniform in $\G(n,d+1)$. 
 
 Finally, we prove that $G_d\subseteq G_{d+1}$ with probability at least $1-5\eta$ for all sufficiently large $n$. We use $G_{d+1}\hookrightarrow \bar G$ and $G_{d+1}\hookrightarrow H$ to denote the events that the coupling procedure sets $G_{d+1}$ to be $\bar G$, and $H$,  respectively. We use $G_{d+1}\hookrightarrow \B$ to denote the event that $\bar G\in \G(n,d+1)\setminus\B$ but $G_{d+1}$ is set to be some graph $G''\in \B$. Note that $G_d\subseteq G_{d+1}$ if $G_{d+1}\hookrightarrow \bar G$. Thus, it is sufficient to show that the probability that $G_{d+1}\hookrightarrow H$ or $G_{d+1}\hookrightarrow  \B$ is at most $5\eta$.
  
  First we see that
 \begin{align*}
 \pr(G_{d+1}\hookrightarrow H\wedge \bar G\in \B)&=\sum_{G\in\G(n,d)}\sum_{G':\substack{G'\sim G\\ G'\in\B}}\frac{\sigma_d}{Z(G)}\left(1-(1-\eta)\frac{Z(G)}{\overline{Z}}\right)\\
 &\hspace{-1cm}\le \sum_{G\in\G(n,d)}\sum_{G':G'\sim G}\frac{\sigma_d}{Z(G)}\left(1-(1-\eta)\left(1-3C\frac{d^2}{n}\right)\right) \quad \mbox{(by Lemma~\ref{lem:bounds})}\\
 &\hspace{-1cm}\le 2\eta\sum_{G\in\G(n,d)}\sigma_d\sum_{G':G'\sim G}\frac{1}{Z(G)}=2\eta,
 \end{align*}
 where the last inequality holds by the definition of $\eta$, the assumptions that $\eta=o(1)$ and $n$ is sufficiently large.
 Similarly,
 \begin{align*}
&\hspace{-3.2cm} \pr\Big(\big(G_{d+1}\hookrightarrow H\ \mbox{or}\ G_{d+1}\hookrightarrow \B\big)\wedge \bar G\notin \B\Big) \\
&\hspace{0.5cm} =\sum_{G\in\G(n,d)} \sum_{\substack{G': G'\sim G\\ G'\notin\B}} \frac{\sigma_d}{Z(G)}\left(1-(1-\eta)\frac{Z(G)}{\overline{Z}} \cdot \frac{\underline{Y}}{Y(G')}\right)  \\
&\hspace{0.5cm} \le 3\eta \sum_{G\in\G(n,d)} \sum_{\substack{G': G'\sim G\\ G'\notin\B\cup\B'}} \frac{\sigma_d}{Z(G)}  +\sum_{G\in\G(n,d)} \sum_{\substack{G': G'\sim G\\ G'\in \B'}} \frac{\sigma_d}{Z(G)},
 \end{align*}
 as for every $G'\notin \B\cup\B'$, $1-(1-\eta)\frac{Z(G)}{\overline{Z}} \cdot \frac{\underline{Y}}{Y(G')} \le \eta + 3\alpha + 6C d^3/n \le 3\eta$ by Lemma~\ref{lem:bounds}, the definition of $\eta$, the assumptions that $\eta=o(1)$ and $n$ is sufficiently large; and for $G'\in \B'$ we use the trivial upper bound $1-(1-\eta)\frac{Z(G)}{\overline{Z}} \cdot \frac{\underline{Y}}{Y(G')}\le 1$.
 Since $|\{G': G'\sim G, G'\notin\B\cup\B'\}|\le |\{G': G'\sim G\}|=Z(G)$, and $\sigma_d\cdot |\G(n,d)|=1$, the first double summation above is at most $3\eta$. By~\eqn{Z-conc}, the second double summation above is equal to
 \begin{align*}
& \left(1+O(d^2/n)\right)\frac{\sigma_d}{Z^\star} |\{(G,G')\in \G(n,d)\times \B': G\sim G'\}|= \left(1+O(d^2/n)\right)\frac{\sigma_d}{Z^\star} d^-(\B')\\
& \hspace{0.5cm} \le \left(1+O(d^2/n)\right)  \frac{\sigma_d}{Z^\star} \eta D \quad \mbox{(by Lemma~\ref{lem:degree-B})} \\
& \hspace{0.5cm}=  \left(1+O(d^2/n)\right)  \frac{\sigma_d}{Z^\star} \eta |\G(n,d)| \ex Z \le 2\eta,
 \end{align*}
 where the last equality above holds because $Z^\star\sim \ex Z$ and $\sigma_d |\G(n,d)|=1$. Combining all cases above, we know that the probability that $G_{d+1}\hookrightarrow \bar G$ (and thus $G_d\subseteq G_{d+1}$) is at least $1-5\eta$.\qed
 
 \subsection{Proof of Theorem~\ref{thm:sandwich2}}
 
 Suppose $d\to\infty$ and $d=o(n^{1/3})$. Then there exists $\alpha=o(1)$ such that $\eta(\alpha)=o(1)$. Theorem~\ref{thm:sandwich2} follows by Theorem~\ref{thm:sandwich3} with such a choice of $\alpha$.  \qed
 
 \subsection{Proof of Theorem~\ref{thm:sandwich}. } 
 Suppose $d\to\infty$, $d=O(n^{1/7}/\log n)$ and $n$ is sufficiently large.
 Let $\alpha=1/d$. It follows now that 
 $\eta=O(1/d)$.  We prove that for each $1\le j\le \lfloor \eps_n d\rfloor$, there is a coupling $(G_d,\ldots, G_{d+j})$ where $G_i\sim \G(n,i)$ for every $d\le i\le d+j$ and with probability at least $1-5j\eta$, $G_d\subseteq G_{d+1}\subseteq\cdots\subseteq G_{d+j}$. Then Theorem~\ref{thm:sandwich} follows by taking $j=\lfloor\eps_n d\rfloor$.
 
 We prove by induction on $j$. The base case $j=1$ follows directly
  by Theorem~\ref{thm:sandwich3} with our choice of $\alpha$.
  Suppose the statement holds for some $1\le j<\lfloor \eps_n d\rfloor$. Let $\pi_j$ be the joint probability distribution of such a coupling $(G_d,\ldots, G_{d+j})$. Again, by
  Theorem~\ref{thm:sandwich3},
 there is a coupling $(G_{d+j}, G_{d+j+1})$ where $G_{d+j}\sim \G(n,d+j)$, $G_{d+j+1}\sim \G(n,d+j+1)$ and with probability at least $1-5\eta$, $G_{d+j}\subseteq G_{d+j+1}$. Let $\pi$ denote the joint probability of this coupling $(G_{d+j}, G_{d+j+1})$. We construct a coupling $(G_d,\ldots, G_{d+j+1})$ by first sampling $(G_d,\ldots, G_{d+j})$ according to the distribution $\pi_j$, and then sampling $G_{d+j+1}$ according to the conditional probability $\pi(G_{d+j+1}\mid G_{d+j})$. The resulting coupling satisfies the required marginal distribution conditions. Moreover, the probability that either $G_d\subseteq\cdots\subseteq G_{d+j}$ fails or $G_{d+j}\subseteq G_{d+j+1}$ fails is at most $5j\eta+5\eta=5(j+1)\eta$ by  the union bound. The assertion follows by induction.\qed
 
\section{Future research} 
 As mentioned in Remark~\ref{remark}, the error $d^3/n$ in Corollary~\ref{cor:subgraph} can be improved to $d^2/n$ if we apply~\cite[Theorem 6]{gao2020triangles} and go through more involved calculations. Another approach is to improve the error in Theorem~\ref{thm:subgraph}, which is of independent interest and can lead to improvements of many other existing results on subgraphs of $\G(n,d)$. 
 
 We solved Conjecture~\ref{con:sandwich} for a certain range of $d_1$ by simultaneously coupling a sequence of random regular graphs. This is certainly not necessary, and is the cause of the restrictions on $d_1$ in Theorem~\ref{thm:sandwich0}. A more direct approach would be to prove concentration of the number of $k$-factors in $\G(n,d)$. This would significantly relax the restrictions on $d_1$, and is itself of independent interest.
 
 \section*{Acknowledgement}
 
 We thank several anonymous referees who carefully read the manuscript and gave many helpful suggestions, which lead to a significant improvement in the presentation.


\section*{Appendix}

\begin{itemize}
\item Justify~\eqn{Exy}
\end{itemize}
 \includegraphics[scale=0.6]{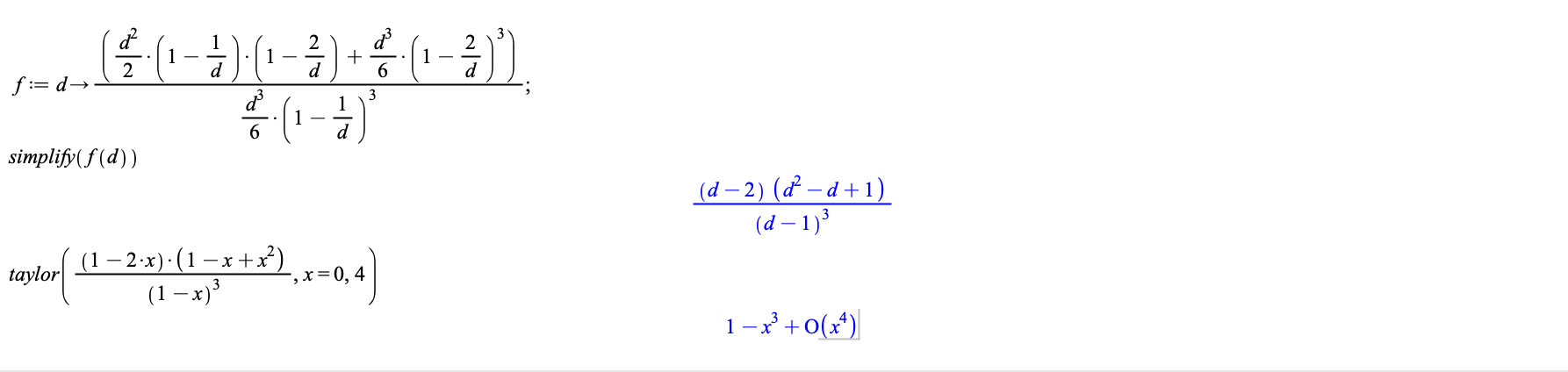}  

\begin{itemize}
\item Justify~\eqn{phi-value}
\end{itemize}
 \includegraphics[scale=0.5]{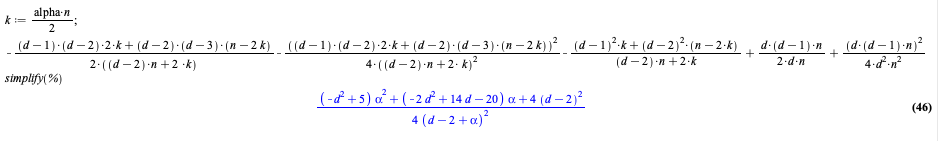}

\begin{itemize}
\item Justify~\eqn{asymp-equation}
\end{itemize}
We verify
\begin{equation}
\sqrt{\frac{n}{e\bar\delta (\bar\alpha n/2)(n-\bar\alpha n)}}  \rho_2(n,d,\bar\alpha)=\left(1+\frac{1}{6d^3}+O(\xi)\right)\rho_1(n,d)^2(\bar\alpha n/e)^{\bar\alpha n/2}. \lab{eq}
\end{equation}
Recall that $\xi=d^{-4}+\frac{d^3}{n}+\sqrt{\frac{d}{n}}\log^3 n$ and $\bar\alpha=1/d$ and
\[
\rho_1(n,d)=\left(\frac{e}{n}\right)^{n/2} \left(\frac{d-1}{d}\right)^{(\frac{d-1}{2})n}d^{\frac{n}{2}}\exp\left(\frac{1}{4}\right) ,
\]
and
\[
\rho_2(n,d,\bar\alpha)=\left(\frac{e}{n}\right)^{(1-\frac{\alpha}{2})n} \left(\frac{d-2+\alpha}{d}\right)^{(\frac{d-2+\alpha}{2})n}d^{\frac{\alpha}{2}n}\left(d-1\right)^{(1-\alpha)n}\exp\left(\phi(d,\bar\alpha)+O(n^{-1})\right).
\]
it is easy to check that all exponential terms cancel exactly from both sides of~\eqn{eq}. By Corollary~\ref{cor:subgraph} with $\alpha=\bar\alpha$ (See Maple calculations and expansions below), 
\[
\phi(d,\bar\alpha)= \frac{4d^2-10d+5}{4(d-1)^2}=1-\frac{1}{2d}-\frac{3}{4d^2}-\frac{1}{d^3}+O(d^{-4}).
\]
 \includegraphics[scale=0.6]{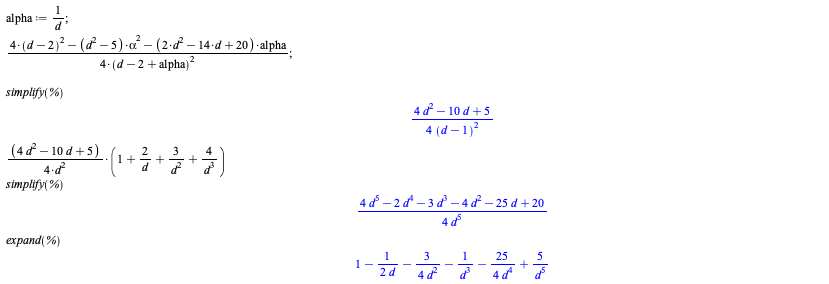}

The polynomially bounded term on the left hand side of~\eqn{eq} is
\begin{align*}
&\sqrt{\frac{n}{e\bar\delta (\bar\alpha n/2)(n-\bar\alpha n)}} \exp\left(\phi(d,\bar\alpha)\right)=\sqrt{\frac{d-1}{e(d-2)}}\exp\left(1-\frac{1}{2d}-\frac{3}{4d^2}-\frac{1}{d^3}+O(d^{-4})\right)\\
&\hspace{0.5cm}=\exp(-1/2) \exp\left(\frac{1}{2d}+\frac{3}{4d^2}+\frac{7}{6d^3}+O(d^{-4})\right) \exp\left(1-\frac{1}{2d}-\frac{3}{4d^2}-\frac{1}{d^3}+O(d^{-4})\right)\\
&\hspace{0.5cm}=\exp\left(\frac{1}{2}+\frac{1}{6d^3}+O(\xi)\right).
\end{align*}
See Maple expansion of $\sqrt{(d-1)/(d-2)}$ below where $x=1/d$:

\noindent \includegraphics[scale=0.6]{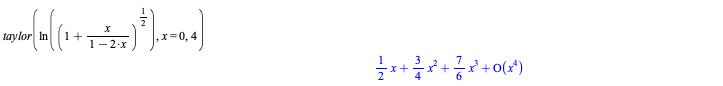}  

\noindent The polynomially bounded term on the right hand side of~\eqn{eq} is
\[
\left(1+\frac{1}{6d^3}+O(\xi)\right) \exp(1/2)=\exp\left(\frac{1}{2}+\frac{1}{6d^3}+O(\xi)\right).
\]

 \end{document}